\documentclass{aptpub}

\authornames{Ora E. Percus and Jerome K. Percus} 
\shorttitle{Maximum of next neighbor walk} 

\numberwithin{equation}{section}  

\newcommand{\lf}{\left}
\newcommand{\rt}{\right}

\newcommand{\be}{\begin{equation}}
\newcommand{\ee}{\end{equation}}

\newcommand{\al}{\alpha}

\newcommand{\de}{\delta}

\newcommand{\la}{\lambda}

\newcommand{\tht}{\theta}

\newcommand{\coef}{\text{coef}}
\newcommand{\sech}{\text{sech}}
\newcommand{\Max}{\text{Max}}
\newcommand{\Var}{\text{Var}}

\newcommand{\half}{\frac{1}{2}}

\newcommand{\tf}{\tilde{f}}
\newcommand{\tg}{\tilde{g}}

\newcommand{\tP}{\tilde{P}}

\newcommand{\tit}{\textit}

\begin{document}

\title{The Maximum of a Symmetric Next Neighbor Walk on the Non-Negative Integers} 

\authorone[New York University]{Ora E. Percus} 

\addressone{Courant Institute, 251 Mercer Street, New York, NY  10012} 

\authortwo[New York University]{Jerome K. Percus} 

\addresstwo{Courant Institute and Physics Department, 251 Mercer Street, New York, NY  10012} 

\begin{abstract}

We consider a one-dimensional discrete  symmetric random walk with a 
reflecting boundary at the origin.  Generating functions are found for 
the 2-dimensional probability distribution $P\{S_n=x, \max_{1\le j \le n} 
S_n=a\}$ of being at position $x$ after $n$ steps, while the maximal 
location that the walker has achieved during these $n$ steps is $a$.  
We also obtain the familiar (marginal) 1-dimensional distribution for 
$S_n=x$, but more importantly that for $\max_{1\le j\le n}\; S_j=a$ 
asymptotically at fixed $a^2/n$.  We are able to compute and compare 
the expectations and variances of the two one-dimensional distributions, 
finding that they have qualitatively similar forms, but differ 
quantitatively in the anticipated fashion.

\end{abstract}



\section{Introduction}

Non-Markovian chains constitute a field of increasing activity.  A dominant 
philosophical motif is that of a hidden Markovian chain~\cite{BP66}, a marginal 
process on a higher dimensional state space.  The analysis of sequences in 
biopolymers \cite{Pe02} as  hidden Markov chains is a primitive version with a small 
underlying state space. 

We were led to consider the problem analyzed in this paper during a study of reinforced 
random walks, next neighbor on a one-dimensional half lattice.  The aim of this paper is to 
find the distribution of $A_n\equiv \max_{1\le i\le n}\, S_i$,  
where $S_i$ is the location of the walker after $i$ steps.
There are numerous ways to  
solve this problem, but we intentionally want to choose one that is extendable to a class 
of reinforced random walks, namely the hidden Markov viewpoint mentioned above.  Before doing 
so, however, it is worth asking what sort of qualitative behavior to expect.   We of course 
will have, asymptotically in $n$, $E(S_n) \propto n^{1/2}$, but the maximum 
sojourn after $n$ steps,  must exceed or equal $S_n$.  How much more? 
But $\{A_i\}$ rectifies the fluctuation in $\{S_i\}$, and hence $A_n$ might be
expected to have a variance, highly reduced from that of $S_n$. How much less?
The limited objective of this paper is to answer these questions by first 
computing the 2-dimensional  distribution of $(S_n, A_n)$ as $n$ varies and 
then the distribution of  the  r.v $A_n$.

\section{Distribution and Moments of $S_n$}

The basic system that we analyze is that of a random walk on the integer lattice 
$x \ge 0$.  The jump $X_i$ at the $i^{th}$ step is next neighbor
\be\label{eq2.1}
X_i = \pm 1
\ee
and the walker starts at the origin, so that its location after $n$ steps is
\be\label{eq2.2}
S_n = \sum^n_{i=1} \, X_i.
\ee
Let us first review the properties of the distribution function
\be\label{eq2.3}
P_n(x) = P(S_n=x).
\ee
We confine our attention to a symmetric walk reflected at the origin, so that
\be\label{eq2.4}
\begin{aligned}
P \Big\{ X_i = \pm 1 \Big| S_{i-1} \neq 0\Big\} &=\half   \\
P \Big\{ X_i =  1 \Big| S_{i-1}=0 \Big\}& =1 . 
\end{aligned}
\ee
The first jump must be from $x=0$ to $x=1$, and so we can take as initial condition
\be\label{eq2.5}
P_1(x) = \de_{x,1}\ .
\ee

The analysis of \eqref{eq2.3} under (\ref{eq2.4}, \ref{eq2.5}) is routine.  We have
\be\label{eq2.6}
\begin{aligned}
P\lf\{ S_1 = x \rt\} &= \de_{x,1} \ ,  \\
P\lf\{ S_n = 0 \rt\} &= \half \: P\lf\{ S_{n-1} =1\rt\} \quad \text{ for }\quad n \ge 2 \ ,  \\
P\lf\{ S_n = 1 \rt\} &= \half \: P\lf\{ S_{n-1} =2\rt\}\;  +\; P \lf\{ S_{n-1} =0 \rt\} \\
P\lf\{ S_n = x \rt\} &= \half \: P\lf\{ S_{n-1} = x+1\rt\}  + \half\: P \lf\{ S_{n-1} =x-1 \rt\} 
\quad \text{for } \quad x \ge 2\ ,
\end{aligned}
\ee
readily solved (index and argument must have the same parity) as
\be\label{eq2.7}
\begin{aligned}
P_{2n} (0) &= \frac{1}{2^{2n}} \binom{2n}{n} \\
P_{2n} (2x) &= \frac{2}{2^{2n}} \binom{2n}{n-x} \quad\text{for }\; x>0\\
P_{2n+1} (2x+1) &= \frac{1}{2^{2n}} \binom{2n+1}{n-x} \quad\text{for }\; x \ge 0.
\end{aligned}
\ee
Observe that \eqref{eq2.7} can also be obtained directly  from a non-reflecting walk from 
the origin to $\pm x$---a trivial combinatorial problem---by reflecting all subwalks on the 
negative axis to the positive axis.  This is because the probability of a walker arriving
at the origin, then jumping to $\pm 1$ is $1$, as in the reflecting case.

Our definition of reflection does not correspond to that of Feller$^{[3]}$, p.~436 and 
Tak\'{a}cs$^{[8]}$ p.~19 where the walker is not allowed to pass a boundary at $x=\half$. 
Instead, when the walker is at $x=1$, the next step takes it to $x=2$ with probability 
$\half$ or it stays at $x=1$ with probability $\half$.  However, Kac$^{[5]}$ and 
Percus$^{[6]}$
treat this walk as a Markov chain with $2\times 2$ transition matrix, equivalent to 
what we do here.

Mean and variance are the leading properties of a random walk, and by direct summation,
one readily finds that
\be\label{eq2.8}
\begin{aligned}
E (S_{2n}) &= \frac{2n}{2^{2n}} \binom{2n}{n} \\
Var (S_{2n}) &= 2n \lf(1- 2n \lf[ \binom{2n}{n}\bigg/2^{2n} \rt]^2\rt),
\end{aligned}
\ee
with a similar result for $S_{2n+1}$.  In both cases, use of the Stirling approximation
shows directly that 
\be\label{eq2.9}
\begin{aligned}
\lim_{n\to\infty} E(S_n)/n^{1/2} &= \sqrt{\frac{2}{\pi}} \\
\lim_{n\to\infty} Var (S_n)/n &= 1-\frac{2}{\pi} \sim  0.36
\end{aligned}
\ee
establishing a standard against which other properties of the walk can be 
compared---the main objective of this paper.

\section{The Joint Distribution $P\{ S_n=x, A_n =a\}$}

Consider then a random walk on the integer lattice $x \ge 0$, $a\ge 1$ with joint 
distribution defined by
\be\label{eq3.1}
\begin{aligned}
&P_n (x,a) \equiv P \lf\{ S_n =x, \;\; A_n=a \rt\}\\
&\text{where } \;S_k = \sum^k_{i=1} \, X_i,\;\;  \qquad X_i = \pm 1.
\end{aligned}
\ee
As in \eqref{eq2.4} we deal with a symmetric random walk reflected at the origin,
and the walk starts at the origin
\be \label{eq3.2}
P_1(x,a) = \de_{x,1} \, \de_{a,1}.
\ee
Since $a$ has not changed from its prior value when $x<a$, we have
\be\label{eq3.3}
\begin{aligned}
P_{n+1} (x,a) &= \half \: P_n (x-1,a) + \half \: P_n(x+1,a) \quad\text{for}\quad
1 < x< a \\
P_{n+1} (1,a) &= P_n (0,a) + \half \: P_n (2,a) (1-\de_{a,1}) \\
P_{n+1} (0,a) &= \half \:P_n(1,a)
\end{aligned}
\ee
But $a$ increases from its prior value with probability $1/2$ when $x=a \ge 1$,
\be\label{eq3.4}
a\ge 2 : P_{n+1} (a,a) = \half \:P_n(a-1, \, a-1) + \half \:P_n(a-1, a).
\ee
We can now combine \eqref{eq3.1}, \eqref{eq3.2}, \eqref{eq3.3}, \eqref{eq3.4} on 
the space defined by 
\be\label{eq3.5}
0 \le x \le a, \qquad\qquad a \ge 2
\ee
obtaining, for $n\ge 1$,
\be\label{eq3.6}
\begin{aligned}
P_{n+1} (x,a) &= \half \lf(1+ \de_{x,1} - \de_{x,a+1}\rt) P_n (x-1,a) \\
&\qquad + \half  \: P_n \, (x+1,a) + \half \: \de_{x,a} \: P_n(a-1, a-1) (1-\de_{a,1})
\end{aligned}
\ee
and initial condition \eqref{eq3.2}.  Note that the condition $P_n(x,a) =0$
for $x>a$, satisfied initially, is automatically satisfied under iteration of 
\eqref{eq3.6}.

Our task is now to solve \eqref{eq3.6}, which we do in standard fashion by first
introducing the generating function, convergent for $|\la| <1$,
\be\label{eq3.7}
\begin{aligned}
P(\la, x, a) &= \sum^\infty_{n=1} \la^n \, P_n (x,a)\\
&= \la \, P_1(x,a) + \sum^\infty_{n=1} \la^{n+1} \, P_{n+1} (x,a).
\end{aligned}
\ee
It follows at once from \eqref{eq3.6} that
\be\label{eq3.8}
\begin{aligned}
P(\la, x,a) &= \la\, P_1 (x,a) + \frac{\la}{2} \lf(1+\de_{x,1} - \de_{x, a+1}\rt) 
P( \la, x-1,a)\\
&\qquad + \frac{\la}{2} \:P (\la, x+1,a) + \frac{\la}{2} \: \de_{x,a}
\: P(\la, a-1, a-1) (1-\de_{a,1}).
\end{aligned}
\ee
Further simplification is then achieved by going over to the double generating 
function 
\be\label{eq3.9}
\tP (\la, u,a) \equiv \sum^\infty_{n=1} \sum^\infty_{x=0} \la^n \, u^x
\, P_n (x,a) = \sum^a_{x=0} P(\la, x,a)\,u^x
\ee
where we have used the fact that $P_n(x,a)=0$ for $x>a$, and this also establishes
that $\tP(\la, u,a)$ is a polynomial in $u$ of degree $a$, thereby convergent for 
all $u$.  Summing \eqref{eq3.8} over $x$, with weight $u^x$, we find after minor
algebra that
\be\label{eq3.10}
\begin{aligned}
&\lf( u^2 - \frac{2u}{\la} +1\rt) \tP(\la,u,a) = -2u^2 \, \de_{a,1} \\
&\qquad \qquad  + \lf(1-u^2\rt) P(\la,0,a) +u^{a+2} \, P(\la, a,a)
- u^{a+1} \, P(\la, a-1, a-1)(1-\de_{a,1}).
\end{aligned}
\ee
Solving \eqref{eq3.10} is fairly straightforward.  First, take the special case
$a=1$:
\be\label{eq3.11}
\lf(u^2-\frac{2u}{\la} +1\rt) \tP(\la,u,1) =-2u^2 + (1-u^2) \, P(\la,0,1)
+ u^3 \, P(\la, 1,1),
\ee
and introduce the zeroes of $u^2-\frac{2u}{\la} +1=0$:
\be\label{eq3.12}
u_1 =\tht= \lf(1-\sqrt{1-\la^2}\rt)\Big/\la, \qquad u_2=\frac{1}{\tht} 
= \lf(1+ \sqrt{1-\la^2}\rt)\Big/ \la.
\ee
Taking $u=\tht$, and then $u=1/\tht$ in \eqref{eq3.11}, we have 
\be\label{eq3.13}
\begin{aligned}
0 &= -2\tht^2 + \lf(1-\tht^2\rt) P(\la, 0,1)+ \tht^3 \, P(\la,1,1)\\
0 &= -2\tht^{-2} + \lf(1-\tht^{-2}\rt) P(\la, 0,1)+ \tht^{-3} \, P(\la,1,1)
\end{aligned}
\ee
and on eliminating $P(\la,0,1)$,
\be\label{eq3.14}
P(\la, 1,1) =2 \lf( \tht + \tht^{-1} \rt) \bigl/ \lf( \tht^2 + \tht^{-2} \rt) 
= \la \bigl/ \lf( 1-\la^2/2\rt).
\ee
The case $a>1$ can be treated the same way.  Using \eqref{eq3.10}
\be\label{eq3.15}
\begin{aligned}
0 &= \lf(1-\tht^2\rt) P(\la, 0,a) + \tht^{a+2} \, P(\la,a,a) - \tht^{a+1}
\, P(\la, a-1, a-1)\\
0 &= \lf(1-\tht^{-2}\rt) P(\la, 0,a) + \tht^{-a-2} \, P(\la,a,a) - \tht^{-a-1}
\, P(\la, a-1, a-1)
\end{aligned}
\ee
and eliminating $P(\la,0,a)$,
\be\label{eq3.16}
P(\la,a,a) =\Bigl[ \Bigl(\tht^a + \tht^{-a}\Bigr)\!\Bigl/\! 
\lf(\tht^{a+1} + \tht^{-(a+1)}\rt)\Bigr]
P(\la, a-1, a-1).
\ee
Starting with \eqref{eq3.14} and iteratively applying \eqref{eq3.16}, we conclude that
\be\label{eq3.17}
P(\la, a,a) = 2 \frac{\tht+ \tht^{-1}}{\tht^{a+1} + \tht^{-(a+1)}}
= \frac{4/\la}{\tht^{a+1} + \tht^{-(a+1)}}
\ee
valid as well for $a=1$, and leading via the first equality 
of \eqref{eq3.14} and the first of \eqref{eq3.17} to (note the 
convention that $P(\la,0,0)=0$)
\be\label{eq3.18}
P(\la, 0,a) = \frac{2\tht^2}{1-\tht^2} \:\de_{a,1} + \frac{\tht^{a+1}}{1-\tht^2}
\:P(\la,a-1, a-1)(1-\de_{a,1}) - \frac{\tht^{a+2}}{1-\tht^2} \:
P(\la, a,a) \ .
\ee 
The net effect, substituting back into \eqref{eq3.10}, is that 
\be\label{eq3.19}
\begin{aligned}
&\tP (\la, u,a) \lf(u^2 -\frac{2u}{\la}+1\rt) = \\
& \qquad \lf( \lf(1-u^2\rt) \frac{\tht^{a+1}}{1-\tht^2} - u^{a+1}\rt) 
\lf\{P (\la, a-1, a-1)(1-\de_{a,1})+ 2\de_{a,1}\rt\}\\
&\hspace{2.0in} - \lf(\lf(1-u^2\rt) \frac{\tht^{a+2}}{1-\tht^2} -u^{a+2}\rt)
P(\la,a,a).
\end{aligned}
\ee
or
\be\label{eq3.20}
\begin{aligned}
(u-\tht) \lf(u-\frac{1}{\tht}\rt) \tP(\la,u,a) &= \lf[ 
\lf(1-u^2\rt) \frac{\tht^{a+1}}{1-\tht^2} - u^{a+1} \rt]
\lf[ 2 \frac{\tht + \tht^{-1}}{\tht^a + \tht^{-a}}\rt]
- \lf[ \lf(1-u^2\rt) \frac{\tht^{a+2}}{1-\tht^2} - u^{a+2} \rt] \\[2mm]
&\qquad \times \lf[ 2 \frac{\tht + \tht^{-1}}{\tht^{a+1} + \tht^{-(a+1)}}\rt]
\end{aligned}
\ee
[From \eqref{eq3.20} we conclude
\begin{align*}
P\{ A_n=a\} &= \text{ coef of }\; \la^n \text{ in }\; \frac{1}{1-\la} 
\lf \{
\frac{2}{\tht^a + \tht^{-a}} - \frac{2}{\tht^{a+1} + \tht^{-(a+1)} }
\rt\}
\intertext{and since}
P\{ A_{n-1}=a-1\} &= \text{ coef of }\; \la^{n-1} \text{ in }\; \frac{1}{1-\la} 
\lf \{
\frac{2}{\tht^{a-1} + \tht^{-(a-1)}} - \frac{2}{\tht^a + \tht^{-a} }
\rt\}
\end{align*}
we also find that
\be\label{eq3.21}
\begin{aligned}
&P \lf\{ \text{ $A_n =a$ \;  for the first time in the $n^{th}$ step } \rt\} \\
&\; = \coef \; \la^n \;\text{ in }\; \frac{1}{1-\la} \lf\{ \frac{2}{\tht^a + \tht^{-a}}
-\frac{2}{\tht^{a+1} + \tht^{-(a+1)}}\rt\}\\
&\; -\coef  \; \la^{n-1}\; \text{ in }\; \frac{1}{1-\la} 
\lf\{ \frac{2}{\tht^{a-1} + \tht^{-(a-1)}} - \frac{2}{\tht^a + \tht^{-a}} \rt\}
\end{aligned}
\ee
(see (4.1)--(4.3) for details)   \hfill ]

\section{ The Limiting Moments of $A_n$}

Our objective is to examine the characteristics of the Non-Markovian random variable
$A_n \equiv \Max_{1 \le j \le n} \, S_n$, which  of course corresponds to obtaining the marginal
distribution in which $P\lf\{ S_n =x, \;\max_{1\le i \le n} \; S_i=a \rt\}$ is summed
over $x$.  The complementary marginal, summed over $a$, is just the usual Markovian 
walk of $P\lf\{ S_n=x\rt\}$, whose solution was given in Sec.~2.

We have seen in \eqref{eq3.20}  that
\be\label{eq4.1}
\begin{aligned}
&\lf(u^2 - \frac{2u}{\la} +1 \rt) \tP(\la, u,a) 
= \lf( \lf(1-u^2\rt) \frac{\tht^{a+1} }{1-\tht^2} - u^{a+1} \rt) 
\lf( 2 \frac{ \tht+ \tht^{-1} }{ \tht^a + \tht^{-a} } \rt)\\
&\hspace{2in} - 
\lf( \lf(1-u^2\rt) \frac{\tht^{a+2} }{ 1-\tht^2} - u^{a+2} \rt)
\lf(2 \frac{ \tht + \tht^{-1} }{ \tht^{a+1} + \tht^{-(a+1)} }\rt)
\end{aligned}
\ee
The generating function for the marginal distribution of $A_n$ is then found by 
summing over $S_n=x$, equivalent to setting $u=1$ in \eqref{eq4.1}:
\be\label{eq4.2}
\lf(1-\frac{1}{\la}\rt) \tP(\la,1,a) = \frac{\tht + \tht^{-1}}{\tht^{a+1} - \tht^{-(a+1)}}
- \frac{\tht + \tht^{-1}}{\tht^a + \tht^{-a}}
\ee
or, since $\tht + \tht^{-1} = 2/\la$,
\be\label{eq4.3}
(1-\la) \, \tP(\la,1,a) = \frac{2}{\tht^a+ \tht^-a} - \frac{2}{\tht^{a+1} +\tht^{-(a+1)}}.
\ee
It is then simple to construct as well the double generating function
\be\label{eq4.4}
\begin{aligned}
Q(\la, z) 
&\equiv \sum^{\infty,\infty}_{1,1} \la^n \, z^a \, P\lf\{ A_n =a \rt\}  \\
& = \sum^\infty_1  \, z^a \, \tP (\la, 1,a)\\
&= \frac{1}{1-\la} \sum^\infty_1 \, z^a  
\lf( \frac{2}{\tht^a + \tht^{-a}} - \frac{2}{\tht^{a+1} \tht^{-(a+1)}} \rt) \\
&= \frac{1}{1-\la} 
\lf[\sum^\infty_1 \lf(z^a - z^{a-1}\rt) \frac{2}{\tht^a + \tht^{-a}} + \la \rt]\\
&= \frac{\la}{1-\la } - \frac{1-z}{1-\la}  \sum^\infty_1 z^a 
\frac{2}{\tht^a + \tht^{-a}}.
\end{aligned}
\ee
The factorial moments of $\{A_n\}$ are of course obtained by $z$-differentiations of 
$Q(\la, z)$ at $z=1$, or directly in the fashion of \eqref{eq4.4}
(using the familiar recurrence relation of binomial coefficients)
\be\label{eq4.5}
\text{For } k \ge 1 \quad\sum_n \la^n E \lf\{ \binom{A_n}{k} \rt\} = \frac{1}{1-\la} 
\sum^\infty_{a=1}
\binom{a}{k-1} \frac{2}{\tht^a + \tht^{-a}}
\ee
Our objective is to obtain  the asymptotic form of $E\binom{A_n}{k}$ as $n\to\infty$.  Our 
claim is that this has the same form as the $k^{th}$-moment of $S_n,  (E(S^k_n))$, which we 
know is proportional to $n^{\frac{k}{2}}$.  In other words we want to find the constant
$C_k$ in the postulated relation
\be\label{eq4.6}
\lim_{n\to\infty} \frac{E\binom{A_n}{k}}{n^{k/2}} = C_k \ .
\ee
It is not obvious that this limit exists because $\frac{1}{n^{\frac{k}{2}}} \:E\binom{A_n}{k}$
may have persistent oscillations when $n\to\infty$.  Therefore we will  instead use a generalized
limit in the sense of Abel or Cesaro in which a suitable running average is performed before
the limit is taken.  Prototypical is one form of Abel limit theorem which states that:
\be\label{eq4.7}
\text{If }\;\lim_{n\to\infty} \, a_n = A \quad\text{ then } \quad \lim_{\la\to 1-}
(1-\la) \sum^\infty_{n=1} a_n \, \la^n =A
\ee
This is readily proved by decomposing the sum into two sums:
\be\label{eq4.8}
(1-\la) \sum^\infty_{n=1} a_n\,\la^n = (1-\la) \sum^{\lf[ (1-\la)^{-\half} -1 \rt]}_{n=1}
\; a_n \, \la^n + (1-\la) \sum^\infty_{n= \lf[ (1-\la)^{-\half} \rt]} a_n\, \la^n
\ee
and observing that the first term $\to 0$ while the second $\to A$ as 
$\la\to 1-$.  Equation~\eqref{eq4.7} can be generalized using the 
same decomposition of the sum (as in 4.8) to read:
\be\label{eq4.9}
\text{If }\quad \lim_{n\to\infty} \frac{1}{\binom{n}{p}} \: a_n = C \quad\text{ then }\quad
\lim_{\la\to 1-} (1-\la)^{1+p} \sum^\infty_{n=1} \la^{n-p} \, a_n =C.
\ee
Since $\binom{n}{p}\: \frac{1}{n^p}$ tends to  $\frac{1}{p!}$ as $n\to\infty$,
we therefore define
\be\label{eq4.10}
\sideset{}{_{}^{*}}{\lim}_{n\to\infty} 
n^{-p}\, a_n \equiv \lim_{\la\to 1-} (1-\la)^{1+p}
\sum^\infty_{n=1} \,  \lf(\frac{1}{p!}\rt)  \la^{n-p} \,a_n \ .
\ee
If $\lim_{n\to\infty} \, n^{-p} \, a_n$ exists then $\lim^*_{n\to \infty} \,a_n \, n^{-p}$
has the same value.  However $\lim^*$ may exist even when $\lim$ does not.

We now apply \eqref{eq4.10} to \eqref{eq4.5} to obtain the $\lim^*$ version of \eqref{eq4.6}
\be\label{eq4.11}
\begin{aligned}
\sideset{}{_{}^{*}}{\lim}_{n\to\infty} 
\frac{1}{ n^{\frac{k}{2}} } \, E\lf\{ \binom{A_n}{k} \rt\} 
&= \lim_{\la \to 1-} (1-\la)^{1+ \frac{k}{2}} \sum^\infty_{n=1} 
\frac{1}{ \lf(\frac{k}{2}\rt)!} \:\la^{n-\frac{k}{2}} \: 
E \lf\{ \binom{A_n}{k}\rt\}\\
&= \lim_{\la \to 1-} \lf( \frac{1-\la}{\la} \rt)^{\frac{k}{2}} \sum_{a=1}
\binom{a}{k-1} \:\frac{2}{\tht^a + \tht^{-a}} \frac{1}{\lf(\frac{k}{2}\rt)!}
\end{aligned}
\ee
Set $\tht=e^{-t}$; therefore
$$
\frac{2}{\tht^a + \tht^{-a}} = \frac{1}{\cosh at} \ ,\quad \qquad \frac{1-\la}{\la} = 
\cosh t-1
$$
and \eqref{eq4.11} can be written as 
\be\label{eq4.12}
\sideset{}{_{}^{*}}{\lim}_{n\to\infty} 
\frac{1}{n^{\frac{k}{2}}} \, E \lf\{\binom{A_n}{k} \rt\} = 
\lim_{t\to 0+} \lf( \cosh t-1\rt)^{\frac{k}{2}} \;\frac{1}{\lf(\frac{k}{2}\rt)!}
\sum_{a=1} \binom{a}{k-1} \: \frac{1}{\cosh at}
\ee
Since $\frac{ (\cosh t-1)^{\frac{k}{2} }}{ \lf( \frac{t^2}{2}\rt)^{ \frac{k}{2} } } 
\;\; \underset{t \to 0}{\longrightarrow} \;1$, \eqref{eq4.12} becomes for
$k\ge 1$
\be\label{eq4.13}
\begin{aligned}
\sideset{}{_{}^{*}}{\lim}_{n\to\infty} 
\frac{1}{n^{\frac{k}{2}}} E\lf\{\binom{A_n}{k}\rt\} 
&= \lim_{t\to 0+} \:\frac{t^k}{2^{\frac{k}{2}} \lf(\frac{k}{2}\rt)!} 
\sum_{a=1} \binom{a}{k-1} \: \frac{1}{\cosh at}\\
&= \lim_{t-\to 0+} \lf(\frac{t}{\sqrt{2}}\rt)^k \:\frac{1}{\lf(\frac{k}{2}\rt)!}
\sum_{a=k-1} \binom{a}{k-1} \:\frac{1}{\cosh at}
\end{aligned}
\ee
But
\be\label{eq4.14}
\begin{aligned}
&\lim_{t\to 0} \lf(\frac{t}{\sqrt{2}}\rt)^k \;\frac{1}{\lf(\frac{k}{2}\rt)!} \sum_{a=k-1}
\binom{a}{k-1} \;\frac{1}{\cosh at} = \\
&\quad \lim_{t\to 0} \;\frac{t^k \,2^{-\frac{k}{2}}}{\lf(\frac{k}{2}\rt)!} 
\lf( \sum^{\frac{\al}{t} -1}_{a=k-1} + \sum^{\frac{\beta}{t} -1}_{a= \frac{\al}{t}}
+ \sum^\infty_{a = \frac{\beta}{t}} \rt) 
\frac{\binom{a}{k-1}}{\cosh at} \ .
\end{aligned}
\ee
It can verified that the contribution of the 3$^{rd}$ sum is a function of $\beta$ and $t$ which
$\to 0$ as $\beta \to \infty$ for any $t$ and the contribution of the first sum $\to 0$ as 
$\al \to 0$ for any $t$.  Then also,
the contribution of the 2$^{nd}$ sum involves a Riemann sum which converges to 
a Riemann integral
\be\label{eq4.15}
\frac{1}{2^{\frac{k}{2}} \lf(\frac{k}{2}\rt)! (k-1)!} 
\int^\beta_\al \,\frac{b^{k-1}}{\cosh b} \:db 
\underset{\substack{\al \to 0\\ \beta \to\infty}}{\longrightarrow}
\frac{1}{2^{\frac{k}{2}} \lf(\frac{k}{2}\rt)! \, (k-1)!} 
\int^\infty_0 \frac{b^{k-1}}{\cosh b} \; db
\ee
From \eqref{eq4.13} and \eqref{eq4.15} we conclude that
\begin{align}
\sideset{}{_{}^{*}}{\lim}_{n\to\infty} 
 \, n^{-\frac{k}{2}} \: E\lf\{ \binom{A_n}{k} \rt\} 
&= \frac{1}{2^{\frac{k}{2}}  \lf(\frac{k}{2}\rt)! \; (k-1)!}
\int^\infty_0 \frac{b^{k-1}}{\cosh b} \; db \label{eq4.16}\\
\intertext{or equally well for $k \ge 1$}
\sideset{}{_{}^{*}}{\lim}_{n\to\infty} 
\, n^{-\frac{k}{2}} \: E\lf\{ A^k_n \rt\} 
&= \frac{k}{2^{\frac{k}{2}}  \lf(\frac{k}{2}\rt)! } \int^\infty_0
\frac{b^{k-1}}{\cosh b} \: db. \label{eq4.17}
\end{align}
The most important examples will be (see [2] for evaluation of the 
integrals)
\be\label{eq4.18}
\begin{aligned}
&k=1 \quad 
\sideset{}{_{}^{*}}{\lim}_{n\to\infty} 
n^{-\half} \;E \lf\{A_n\rt\}
= \frac{1}{\sqrt{2}\: \half \sqrt{\pi}}   
\int^\infty_0 \frac{db}{\cosh b}
= \sqrt{\frac{\pi}{2}}\\
&k=2 \quad 
\sideset{}{_{}^{*}}{\lim}_{n\to\infty} 
\frac{1}{n} \;E \lf\{A^2_n\rt\} =
\int^\infty_0 \frac{bdb}{\cosh b}
= 2G = 1.83193 \\
&\text{where $G$ is the Catalan constant (see [4])}
\end{aligned}
\ee
from which
\be\label{eq4.19}
\sideset{}{_{}^{*}}{\lim}_{n\to\infty} 
\, \frac{1}{n} \,\Var \, (A_n)= 0.26113 \dots \ .
\ee
Including the results of \eqref{eq2.9} we have
\begin{alignat*}{2}
&\lim \, \frac{1}{\sqrt{n}} \: E(S_n) \equiv \sqrt{\frac{2}{\pi}} = 0.7979
&\qquad 
&\sideset{}{_{}^{*}}{\lim} 
\frac{1}{\sqrt{n}} \: E(A_n) = \sqrt{\frac{\pi}{2}} = 1.2533\\ 
&\lim \, \frac{1}{\sqrt{n}} \: \Var(S_n) = 0.36
&\qquad 
&\sideset{}{_{}^{*}}{\lim} 
\, \frac{1}{n} \: \Var(A_n) = 0.26113 
\end{alignat*}
From the definitions of $A_n$ and $S_n$ the mean of $A_n$ might have been much 
larger than the mean of $S_n$:  It is not.  The variance of $A_n$ might have been 
much smaller than the variance of $S_n$:  It is not.

A heuristic version of the process used in obtaining \eqref{eq4.16} can be carried
out as well for \eqref{eq4.4} resulting in 
\begin{align*}
Q(\la,z) = \frac{1}{1-\la} \lf(\la + \half \lf(1-z\rt)\rt) 
&- \frac{1-z}{1-\la} \;\frac{1}{2\sech^{-1}\,\la}
\lf\{ \lf(
\psi \lf(\frac{1}{4} \;\frac{\log z}{\sech^{-1} \,\la}\rt) + \frac{3}{4}
\rt) \rt. \\
&\lf.   - \psi \lf(\frac{1}{4} \;\frac{\log z}{ \sech^{-1} \,\la} 
+ \frac{1}{4} \rt) \rt\}
\end{align*}
where $\psi$ is the dilogarithm function.

\section{The limiting distribution of $A_n$}

The moments supply crucial information as to the nature of the distribution
of the random variation $A_n$.  But can we find this distribution 
$$
Q_n(a) \equiv P\{A_n=a\}
$$
in an explicit---and usable---form?  We have seen, and used, the generating
function relation
\be
(1-\la) \:\tP (\la,1,a) = \frac{2}{\tht^a+\tht^{-a}} 
- \frac{2}{\tht^{a+1} + \tht^{-(a+1)}}
\tag{4.3}
\ee
with the consequence that
\be\label{eq5.1}
Q_n(a) - Q_{n-1} (a) = \coef \,\la^n \;\text{ in }\; 
\frac{2}{\tht^a + \tht^{-a}} - \frac{2}{\tht^{a+1} + \tht^{-(a+1)}}
\ee
Our first task will be to find $\coef \, \la^n$ in $\frac{1}{\tht^a+\tht^{-a}}$. 
It is easy to see that 
$$
\frac{2}{\tht^a +\tht^{-a}} = \frac{1}{T_a \lf(\frac{1}{\la}\rt)}
$$
where $T_a$ is the $a^{th}$ Chebyshev polynomial, but this is not very helpful.
However, a simple partial fraction decomposition is completely effective.  We have
\be\label{eq5.2}
\begin{aligned}
\frac{1}{\tht^{-a} + \tht^a} &= \frac{\tht^a}{\tht^{2a} +1}  \\
&= \sum^{2a}_{j=1} \frac{\tht^a_j}{2a \,\tht^{2a-1}_j}
\lf(\frac{1}{\tht- \tht_j}\rt)\\
&= \frac{1}{2a} \;\sum^{2a}_{j=1} \frac{\tht^{a+1}_j }{\tht-\tht_j}\\
\text{where $\tht_j$} &= e^{(i\pi/2a)(2j-1)}\ .
\end{aligned}
\ee
Replacing $\tht_j$ by $1/\tht_j$ does not change the set $\{\tht_j\}$, and so we 
can replace \eqref{eq5.2} by its average over the two forms:
\be\label{eq5.3}
\begin{aligned}
\frac{1}{\tht^{-a} + \tht^a}
&= -\frac{1}{4a}\, \sum^{2a}_1 \frac{\tht^{a+1}_j}{\tht-\tht_j}
+ \frac{\tht^{-(a+1)}_j}{\tht-\tht^{-1}_j} , \\
&= \frac{1}{4a} \sum^{2a}_1 
\frac{\tht^a_j \lf(\tht^{-1}_j - \tht_j\rt) }{ 
\tht + \tht^{-1} - \lf(\tht_j + \tht^{-1}_j\rt)}\\
&= \frac{i\la}{8a} \sum^{2a}_1
\frac{(-1)^j \lf(\tht^{-1}_j - \tht_j\rt)}{
1-\frac{\la}{2}  \lf(\tht^{-1}_j + \tht_j\rt) }
\end{aligned}
\ee
where we have used $\tht_j=\exp (i\pi/2a) (2j-1)$, $\tht^a_j=-i (-1)^j$, and 
$\tht + \tht^{-1}=2/\la$.  It follows at once that 
\be\label{eq5.4}
\text{coef  of $\la^n$ in }\;\frac{1}{\tht^a + \tht^{-a}} = \frac{1}{4a}
\sum^{2a-1}_0 (-1)^j \; \sin \;\frac{\pi}{2a} \; (2j+1) \, \cos^{N-1}
\; \frac{\pi}{2a} \: (2j+1).
\ee

Note that the summation index replacement $j\to 2a-1-j$ in \eqref{eq5.4}
leaves every term unchanged.  Thus, we can replace the summation range by 
its lower half and multiply by 2.
\be\label{eq5.5}
\text{coef $\la^n$ in }\; \frac{1}{\tht^a +\tht^{-a}} = \frac{1}{2a}
\sum^{a-1}_{j=0} (-1)^j \, \sin \: \frac{\pi}{2a} \:(2j +1) \, \cos^{n-1}
\:\frac{\pi}{2a} \; (2j+1) \ .
\ee
But then, the replacement $j\to a-1-j$ in \eqref{eq5.5} multiplies each
term by $(-1)^{n+a}$, with four consequences:
\begin{enumerate}
\item[i)] $\coef \la^n \text{ in } \frac{1}{\tht^a + \tht^{-a}} \neq 0$ only 
if $a\equiv n \pmod{ 2}$
\item[ii)] If $a\equiv n \pmod{2}$, then \eqref{eq5.5} can be reduced to its 
lower half-range (for odd $a$, the summand vanishes at both $(a-1)/2$ and 
$(a+1)/2$):
\be\label{eq5.6}
\coef \;\la^n \text{ in }\; \frac{1}{\tht^a + \tht^{-a}} = \frac{1}{a} 
\sum^{[a/2-1]}_{j=0} (-1)^j \,\sin \:\frac{\pi}{2a}\: (2j+1)\, \cos^{n-1}
\: \frac{\pi}{2a} \: (2j+1)
\ee
\item[iii)] The sum \eqref{eq5.6} is strictly alternating in sign, since
$0 < \frac{\pi}{2a} \:(2j+1) \le \frac{\pi}{2}$\ .
\item[iv)] It also follows from i) that
\be\label{eq5.7}
\begin{aligned}
Q_n(a) - Q_{n-1} (a) &= \de_{a,1} \, \de_{n,1} + 2(-1)^{n+a} \; \coef\;\;
\la^n \text{ in } \;\frac{1}{\tht^b + \tht^{-b}} \\
\text{where } b &= 
\begin{cases}
a &\text{for  $a \equiv n\!\!\pmod{2}$} \\
a+1 &\text{for $a \equiv (n+1)\!\! \mod{2}$}
\end{cases}
\end{aligned}
\ee
\end{enumerate}
We can apply \eqref{eq5.7} at once to \eqref{eq5.1} by making use of the fact
that $\lim_{a\to\infty} \,P\lf\{A_n=a \rt\} \linebreak
=0$.  It then follows from \eqref{eq5.6}
on summing over $n$ that
\be\label{eq5.8}
\begin{aligned}
&\; Q_N(a) = -\sum^\infty_{n=N+1} \lf[Q_n(a) - Q_{n-1} (a)\rt]=\\
&\begin{cases}
\frac{1}{a+1} \sum_j  (-1)^j\;
\frac{\cos^N \frac{\pi}{2(a+1)}\:(2j+1)}{\sin \frac{\pi}{2(a+1)}\: (2j+1)} 
- \frac{1}{a} \sum_j (-1)^j \:
\frac{\cos^{N+1} \frac{\pi}{2a}\:(2j+1)}{\sin\frac{\pi}{a}\:(2j+1)}
&\text{if  $N\equiv a\!\! \mod{2}$}\\
\frac{1}{a+1} \sum_j (-1)^j\; \frac{\cos^{N+1}\frac{\pi}{2(a+1)}\:(2j+1)}{
\sin \frac{\pi}{2(a+1)} \:(2j+1)} - \frac{1}{a} \sum_j (-1)^j
\: \frac{\cos^N \frac{\pi}{2a} \:(2j+1)}{\sin \frac{\pi}{2a}\:(2j+1)}
&\text{if  $N\equiv (a +1)\!\! \mod{2}$}
\end{cases}
\end{aligned}
\ee
which although rather complex has the necessary property of vanishing
when $N<a$.
Note that the sum over $n$ required to obtain \eqref{eq5.8} starts at $N+1$
or $N+2$ depending upon the relative parity of $N$ and $a$, and goes up in 
steps of $2$.

We'll find the limit of a $Q_N(a)$ as $a$ and $N\to\infty$ at  fixed $\gamma$
where
\be\label{eq5.9}
\frac{a^2}{N} = \gamma \:\frac{\pi}{2} + O\lf(\frac{1}{a}\rt).
\ee
Consider the case $N \equiv a \mod{2}$ in \eqref{eq5.8} (the case $N\equiv (a+1)
\mod{2}$ proceeds similarly).
After certain amount of algebra one finds 
\be\label{eq5.10}
\lim_{\substack{ a, N\to\infty\phantom{xxxx} \\
\frac{a^2}{N} = \gamma\: \frac{\pi}{2} + O\lf(\frac{1}{N}\rt)} }
\lf\{ 
\frac{a}{a+1} \:\frac{\cos^N \frac{\pi}{2(a+1)} \:(2j+1)}{ \sin \frac{\pi}{2(a+1)}
\: (2j+1)} - \frac{\cos^{N+1} \frac{\pi}{2a} \:(2j+1)}{\sin \frac{\pi}{2a}
\: (2j+1)} \rt\} =\frac{(2j+1)}{\gamma} \: e^{-\frac{\pi}{4\gamma} \:(2j+1)^2}
\ee

It therefore follows that for fixed $\gamma$
\be\label{eq5.11}
\lim_{a\to\infty} \,a\,Q_N(a) = \sum^\infty_{j=0} (-1)^j \: \frac{(2j+1)}{\gamma}
\:e^{-\frac{\pi}{4\gamma} \:(2j+1)^2}.
\ee
For $\gamma<1$, \eqref{eq5.11} is an alternating series with a decreasing absolute 
value of the $j^{th}$ term.  The absolute ratio of the $j^{th}$ term to the 
$(j-1)^{st}$ term is given by 
$\frac{2j+1}{2j-1}\; e^{-\frac{2\pi j}{\gamma}}$ and 
\be\label{eq5.12}
\begin{aligned}
&\text{hence }\; (1-\al)\:\frac{1}{\gamma} \:e^{-\frac{\pi}{4\gamma}} 
\le \lim_{a\to\infty} a\,Q_N(a) \le \frac{1}{\gamma} \:e^{-\frac{\pi}{4\gamma}} \\
&\text{where }\:\al= 3e^{-\frac{2\pi}{\gamma}} \le 0.0056, \quad
\text{for }\;\gamma \le 1
\end{aligned}
\ee
For $\gamma >1$ the rapid convergence of the series \eqref{eq5.11} quickly 
deteriorates as does the information supplied by the $1^{st}$ term in the series.
However \eqref{eq5.11} does exist and absolutely converges for all $\gamma$;
it is therefore necessary to replace \eqref{eq5.11} by a more rapidly convergent
representation.  This is supplied by a modification of the familiar Poisson
resummation 
\be\label{eq5.13}
\sum^\infty_{-\infty} g(j) = \sum^\infty_{k=-\infty} \tg (k) 
\quad\text{where}\quad \tg(k) \equiv 
\int^\infty_{-\infty} g(x) \,e^{2\pi ikx}\,dx 
\ee
As a special case define $g(j) \equiv \sin \lf(\frac{\pi}{2}\:j\rt) f(j)$ 
so that $\sum^\infty_{j=-\infty} g(j) = \sum^\infty_{n=-\infty} 
(-1)^n\:f(2n+1)$.
It then follows directly from \eqref{eq5.13} that 
\be\label{eq5.14}
\sum^\infty_{n=-\infty} (-1)^n \, f(2n+1) = \frac{1}{2i} \sum^\infty_{k=-\infty}
(-1)^k \, \tf\lf(\frac{1}{4} \, (2k+1)\rt).
\ee
As an example, we find at once that
\be\label{eq5.15}
\gamma^{3/2} \sum^\infty_{n=-\infty} (-1)^n \lf(n+ \half\rt) e^{-\pi \gamma
\lf(n+\half\rt)^2} = \sum^\infty_{k=-\infty} (-1)^k \lf(k+\half\rt) e^{-\frac{\pi}{\gamma}
\lf(k+\half\rt)^2}.
\ee
Hence \eqref{eq5.11} is equivalent to 
\be\label{eq5.16}
\lim a\,Q_n(a) = \sqrt{\gamma} \sum^\infty_{j=0} (-1)^j \;(2j+1) \, e^{-\frac{\pi\gamma}{4}
\lf(2j+1\rt)^2}.
\ee
Equation \eqref{eq5.16} now converges very rapidly as did \eqref{eq5.11} and we 
similarly conclude that
\be\label{eq5.17}
(1-\al) \,\sqrt{\gamma}\:e^{-\frac{\pi\gamma}{4}} \le \lim_{a\to\infty} a\,Q_n(a)
\le \sqrt{\gamma} \:e^{-\frac{\pi\gamma}{4}}
\ee
for $\al=3e^{-2\pi\gamma}\le 0.0056$ when $\gamma \ge 1$.

The general summation device we have used is not unknown in our particular case;
It stems from the fact that \eqref{eq5.11} is recognized as a derivative of the 
Jacobi theta function, which under the Jacobi imaginary transformation is converted
to \eqref{eq5.16}.

\section{Concluding Remarks}

We conclude [\eqref{eq5.11} and \eqref{eq5.16}] that the asymptotic
$\lf(a\to\infty \text{ at constant }\gamma = \frac{2a^2}{\pi n} + 
O\lf(\frac{1}{a}\rt)\rt)$ value of the pointwise
distribution $a\,Q_n(a)=a\,P(A_n=a)$ of the maximum of our random walk has been 
found over the full range of $\gamma$.  Furthermore, a very simple estimate,
\eqref{eq5.12} and \eqref{eq5.17} was obtained with a uniform maximum relative 
error of $\al=0.0056$.  Coupled with the asymptotic $(n\to\infty)$ Abel-smoothed
moments we have presented in \eqref{eq4.16}, \eqref{eq4.17}, a quite complete 
characterization of this process has become available.

\section*{Acknowledgements}

We gratefully acknowledge the invaluable help of an anonymous referee in clarifying
the presentation of this paper.  Conversations with J.~H.~Spencer and S.R.S.~Varadhan
were helpful as well.

\end{document}